\documentclass{article}
\usepackage{latexsym}
\usepackage{amsfonts}
\usepackage{amssymb}
\usepackage{graphicx}
\usepackage{graphicx}
\usepackage{latexsym,amsmath,amscd,amsthm,amsxtra}

\newtheorem{theorem}{Theorem}[section]

\newtheorem{lemma}[theorem]{Lemma}
\newtheorem{proposition}[theorem]{Proposition}

\newcommand{\g}{\mathfrak{g}}

\newcommand{\ak}{\mathfrak{k}}
\newcommand{\ap}{\mathfrak{p}}
\newcommand{\ba}{\mathfrak{a}}
\newcommand{\ab}{\mathfrak{b}}
\newcommand{\al}{\mathfrak{l}}
\newcommand{\as}{\mathfrak{s}}
\newcommand{\am}{\mathfrak{m}}

\newcommand{\AC}{\mathbb{C}}

\newcommand{\D}{\displaystyle}

\newcommand{\pf}{\noindent{\bf Proof.}\ }

\begin{document}

\title{The Weyl Integration Model for KAK decomposition of
Reductive Lie Group}

\date{ }

\author{HU Zhi-guang$^{1}$ $\quad$ YAN Kui-hua$^{1,\,2}$\\
 (1. \it School of Mathematical Sciences, Peking University,
Beijing,
100871; \\  2. \it School of Mathematics and Physics, Zhejiang Normal University, \\
\it Zhejiang Jinhua, 321004)\\
\it(Email: \;zhiguang\_hu@eyou.com;\;
yankh@zjnu.cn) \\
}


\maketitle

\begin{abstract}

\vskip 0.1cm The Weyl integration model presented by An and Wang
can be effectively used to reduce the integration over $G$-space.
In this paper, we construct an especial Weyl integration model for
KAK decomposition of Reductive Lie Group and obtain an integration
formula which implies that the integration of $L^1$-integrable
function over reductive Lie group $G$ can be carried out by first
integrating over each conjugacy class and then integrating over
the set of conjugacy classes.



\vskip 0.1cm {\bf CLC number:} O152.5; \, O186.1

\end{abstract}

\vskip 0.5cm
\section{Introduction}

\vskip 0.3cm \quad\, For a compact connected Lie group $G$, the
classical Weyl integration formula (e.g. ref. Knapp\cite{Kn})
indicates that an integration of any continuous function over $G$
can be reduced to one over its maximal torus. That is to say the
integration can be carried out by first integrating over each
conjugacy class and then integrating over the set of conjugacy
classes. Then for a general Lie group, how to reduce the
integration of functions over it\,?

\vskip 0.3cm In Helgason\cite{He2}, the author presents some
integral formulas related to the Cartan, Iwasawa and Bruhat
decompositions for semisimple Lie groups. In this paper, we will
use the Weyl integration model, which is first introduced by An
and Wang in \cite{AW} and can be used to generalize the reducing
integration idea to integrations over G-spaces, to obtain an
integral formula for KAK decomposition of reductive Lie groups. It
is worthy of indicating that this reducible idea of integrations
can be effectively used to calculate the eigenvalue distribution
of random matrices in random matrix ensemble theory. But here we
do not plan to study those applications in random matrix theory.

\vskip 0.3cm Let $G$ be a Lie group, $X$ a $G$-space, $Y$ a
imbedding submanifold of $X$ and $\sigma$ the $G$-action on $X$.
Suppose that $dx,\,dy$ and $d\mu$ are proper invariant measures
over $X,\,Y$ and $G/K$ respectively. The Weyl integration model
(see section 2 for its details) is a system
$\big(G,\sigma,(X,dx),(Y,dy),(G/K,d\mu)\big)$ satisfying some
proper conditions (see \eqref{orbit-intersect}-\eqref{dimK-Ky}).
Indeed, for a Weyl integration model, the following formula holds
for all $f\in L^1(X,dx)$,
\begin{gather}\label{F:Weylbianhuan}
\int_X f(x)
dx=\frac{1}{d}\int_Y\Big(\int_{G/K}f(\sigma(g,y))\,d\mu([g])\Big)
J(y)\,dy ,
\end{gather}
where $J\in C^\infty(Y)$ and $d$ is the multiplicity of the model.
That is to say that the integration over a $G$-space $X$ can be
converted to the integration by first integrating over each orbit
and then integrating over the orbits space $Y$. Moreover, under
some orthogonal conditions, $J(y)$ can be effectively calculated,
i.e. there is a constant $C$ such that
\begin{gather}\label{P:det2}
J(y)=C\big|\det\Psi_y\big|, \quad \text{for all}\; y\in Y',
\end{gather}
where $\Psi_{y}$ is a mapping induced by G-action $\sigma$ (see
\eqref{Psi}).

\vskip 0.3cm In this paper, we will construct a Weyl integration
model $\Big(K\times K, \sigma, (G,dg), (A, da), \big((K\times
K)\big/M, d\mu\big)\Big)$. The meanings of its members are as
follow. $G$ is a reductive Lie group with Lie algebra $\g$. Let
$\theta$ be the Cartan involution on $\g$. $\g=\ak\oplus\ap$ is
the Cartan decomposition of $\g$, where $\ak$ and $\ap$ are the
$+1$ and $-1$ eigenspaces of $\theta$ respectively. Let $K$ be the
associated maximal compact subgroup of $G$ with Lie algebra $\ak$.

Let $\ba$ be a maximal abelian subspace of $\ap$. Set
$A=\exp(\ba)$, then $A$ is a closed subgroup of $G$. It is known
that $G=KAK$ in the sense that every element in $G$ has a
decomposition as $k_{_1}ak_{_2}$ with $k_{_1},k_{_2}\in K$ and
$a\in A$.

To consider the following group $K\times K$-action
\begin{gather}\label{G-action-varphi}
\sigma:(K\times K)\times A\rightarrow G,\quad
\big((k_{_1},k_{_2}), a\big)\mapsto k_{_1}ak_{_2}^{-1}.
\end{gather}
Let
\begin{gather}\label{M}
M=\big\{(k_{_1},k_{_2})\in K\times K\,\big|\,k_{_1}ak_{_2}^{-1}=a,
\forall\; a\in A\big\}.
\end{gather}
 Then
the $K\times K$-action $\sigma$ can be naturally reduced to a map
$\varphi$, i.e.
\begin{gather}\label{reduced-map-varphi}
\varphi: (K\times K)\big/M\times A\rightarrow G,\quad
\big([(k_{_1},k_{_2})],a\big)\mapsto k_{_1}ak_{_2}^{-1}.
\end{gather}

Let $dg,da$ be the Haar measures on $G$ and $A$ respectively.
 There has a $K\times K$-invariant measure $d\mu$ on
 $(K\times K)\big/M$ since $K\times K$ is compact.
Then we can prove that $\Big(K\times K, \sigma, (G,dg), (A, da),
\big((K\times K)\big/M, d\mu\big)\Big)$ is a Weyl integration
model with finite multiplicities. Thus by \eqref{F:Weylbianhuan},
we obtain an integral formula (see theorem \ref{main-result1} for
details) which reduces the integration over reductive Lie group
$G$ to one first over $(K\times K)\big/M$ and then over its
subgroup $A$. Moreover, by \eqref{P:det2}, we can calculate that
the according $J(y)$ has the following formula
\begin{gather}\label{noncompactJa}
J(a)=C\prod_{\lambda\in
\Sigma^+}\big|\sinh\big(\lambda(H)\big)\big|^{\beta_{_\lambda}},
\end{gather}
where $a=e^H,\, H\in\ba$, $\beta_{_\lambda}$ is the dimension of
the restricted root space $\g_{_\lambda}$ with respect to the
restricted root $\lambda$ in the positive restricted root system
$\Sigma^+$ and $C$ is a proper constant. By comparison, it is
formally analogous to the integration formula for Cartan
decomposition of noncompact semisimple Lie groups (see
\cite{He2}).

\vskip 0.3cm The structure of the paper is as follow. In section
2, we briefly introduce the Weyl integration model and the
Restricted root system and KAK decomposition for a reductive Lie
group. The main results of the paper are presented at the end of
this section. In section 3, we prove the main results.

\vskip 0.5cm
\section{Preliminaries and Main Results}

\vskip 0.3cm \subsection{Weyl integration model}

\quad\, In this subsection, we briefly introduce the Weyl
integration model. One may refer to \cite{AW} for details.

Let $G$ be a Lie group which acts on a $n$-dimensional smooth
manifold $X$. The action is denoted by $\sigma: G\times
X\rightarrow X$. Let $dx$ be a $G$-invariant admissible measure on
$X$. $Y$ is an imbedding submanifold of $X$. Suppose that there is
an admissible measure $dy$ on $Y$, and $X_0 \subset X$, $Y_0
\subset Y$ are closed zero measure subsets of $X$ and $Y$
respectively. Set $X' = X \setminus X_0$, $Y' = Y \setminus Y_0$
and $K=\big\{g\in G\,\big|\,\sigma(g,y)=y, \forall\; y\in
Y\big\}$. Let $K_x=\big\{g\in G\,\big|\,\sigma(g,x)=x\big\}$ be
the isotropic subgroup associated with $x\in X$ and
$O_y=\big\{\sigma(g,y)\,\big|\,g\in G\big\}$ the orbit of $y\in
Y$. Then $K\subset K_y, \forall\; y\in Y$. In the following text
we suppose that
\begin{gather}\label{orbit-intersect}
X'=\bigcup_{y\in Y'} O_y.
\end{gather}
\begin{gather}\label{transversal-condition}
T_y X=T_y O_y\oplus T_y Y, \quad \forall y\in Y',
\end{gather}
which is called the transversal condition.
\begin{gather}\label{dimK-Ky}
\mathrm{dim}K_y=\mathrm{dim}K, \quad \forall y\in Y'.
\end{gather}

\vskip 0.3cm The G-action $\sigma: G\times X\rightarrow X$ can be
reduced to a map $\varphi: G/K\times Y\rightarrow X$ by
$\varphi([g],y)=\sigma(g,y)$ and furthermore to a map which is
still denoted by $\varphi$, i.e. $\varphi: G/K\times Y'\rightarrow
X'$ by restriction. By the above assumption, $\varphi$ is
surjective. Suppose that there is a $G$-invariant admissible
measure $d\mu$ on $G/K$.

\vskip 0.3cm
\begin{proposition}\label{P:localdiffeo}
Suppose the conditions \eqref{orbit-intersect},
\eqref{transversal-condition} and \eqref{dimK-Ky} hold, then
$\varphi: G/K\times Y'\rightarrow X'$ is a local diffeomorphism.
\end{proposition}

Let $G,X,Y,K,\sigma,\varphi,dx,dy,d\mu$ be the above objects. A
\emph{Weyl integration model} is a system
$\big(G,\sigma,(X,dx),(Y,dy),(G/K,d\mu)\big)$, in which we can
choose $X_0$ such that the map $\varphi: G/K\times Y'\rightarrow
X'$ is a finite-sheeted covering map. The number of sheets of the
covering map is called the \emph{multiplicity} of the model. About
a Weyl integration model, we have the following basic theorem.

\vskip 0.3cm
\begin{theorem}\label{T:Weylbianhuan}
If $\big(G,\sigma,(X,dx),(Y,dy),(G/K,d\mu)\big)$ is a Weyl
integration model with multiplicity $d$, then the formula
\eqref{F:Weylbianhuan} holds for all $f\in C^\infty(X)$ with
$f\geq0$ or $f\in L^1(X,dx)$.
\end{theorem}

Let $\as$ be a linear subspace of the Lie algebra $\g$ of $G$,
such that $\g=\ak\oplus\as$, where $\ak$ is the Lie algebra of
$K$. Then $\as$ can be identified with the tangent space
$T_{[e]}(G/K)$ of $G/K$ at point $[e]$ in a natural way:
$$
T_{[e]}(G/K) \cong \as .
$$

By the transversal condition \eqref{transversal-condition}, the
tangent map $d\varphi:T_{[e]}(G/K) \oplus T_yY \rightarrow T_yX$
of $\varphi$ at point $([e], y)$ ($y\in Y'$) can be regarded as a
linear transformation
$$
d\varphi : \as \oplus T_yY \rightarrow T_yO_y \oplus T_yY .
$$
It can be proved that (see \cite{AW})
\begin{gather}
d\varphi=\left(\begin{array}{cc}\Psi_y&0\\0&id
\end{array}\right),
\end{gather}
where $\Psi_y : \as \rightarrow T_yO_y$ is given by
\begin{gather}\label{Psi}
\Psi_y(\xi)=\frac{d}{dt}\Big|_{t=0} \sigma(\exp t\xi,y), \quad
\forall\; \xi\in \as.
\end{gather}

\vskip 0.5cm Suppose that there is a Riemannian structure on $X$
such that the following orthogonal condition holds
\begin{gather}\label{orthogonal-condition}
T_yY\perp T_yO_y, \quad \text{for all}\; y\in Y'.
\end{gather}
Let $dx$ and  $dy$ are the Riemannian measures on $X$ and $Y$
respectively. Then there is a constant $C$ such that
\begin{gather}\label{P:det}
J(y)=C\big|\det\Psi_y\big|, \quad \text{for all}\; y\in Y'.
\end{gather}


\vskip 0.3cm \subsection{Restricted root system and KAK
decomposition for reductive Lie group}

\quad\, In the sense of Knapp (\cite{Kn}, Sec.7.2), a
\emph{reductive Lie group} is a 4-tuple $(G,K,\theta,B)$
consisting of a Lie group $G$, a compact subgroup $K$ of $G$, a
Lie algebra involution $\theta$ of the Lie algebra $\g$ of $G$ and
a nondegenerate $Ad(G)$-invariant $\theta$-invariant bilinear form
$B$ on $\g$ such that (i) $\g$ is a reductive Lie algebra, (ii)
the decomposition of $\g$ into $+1$ and $-1$ eigenspaces under
$\theta$ is $\g=\ak\oplus \ap$, where $\ak$ is the Lie algebra of
K, (iii) $\ak$ and $\ap$ are orthogonal under $B$ and $B$ is
positive definite on $\ap$ and negative definite on $\ak$, (iv)
multiplication, as a map from $K\times \exp(\ap)$ into $G$, is a
diffeomorphism onto, and (v) every automorphism $Ad(g)$ of
$\g^{\AC}$ is inner for $g\in G$, i.e., is given by some $x\in
Int(\g^{\AC})$. $K$ is called the associated \emph{maximal compact
subgroup}, $\theta$ the \emph{Cartan involution} and $B$ the
\emph{invariant bilinear form}. The decomposition of $\g$
(according $G$) in property (iii) (according (iv)) is called
\emph{(global) Cartan decomposition}.

\vskip 0.3cm Now let $G$ is a reductive Lie group with Lie algebra
$\g$. Let $\ba$ be a maximal abelian subspaces of $\ap$, then
$\D\ap=\bigcup_{k\in K}Ad(k)\ba$. Set $\D
\g_{_\lambda}=\big\{X\in\g\,\big|\,(adH)X=\lambda(H)X,\;\text{for
all}\; H\in \ba\big\}$. A nonzero $\lambda\in \ba^*$ is called a
\emph{restricted root} of $\g$ if $\g_{_\lambda}$ is nonzero.
Accordingly $\g_{_\lambda}$ is called a restricted root space. The
set of restricted roots is denoted by $\Sigma$ and Let $\Sigma^+$
be the set of positive restricted roots.

\vskip 0.3cm Reflections in the restricted roots generate the
\emph{Weyl group} $W(\Sigma)$ of $\Sigma$. Denoted by $N_{K}(\ba)$
and $Z_{K}(\ba)$ the normalizer and centralizer of $\ba$ in $K$
respectively, then the Weyl group $W=N_{K}(\ba)\big/Z_{K}(\ba)$
and the Lie algebras of $N_{K}(\ba)$ and $Z_{K}(\ba)$ are
$\am=Z_{\ak}(\ba)$.

\vskip 0.3cm The restricted root space $\g_{_\lambda}$ satisfies
the following basic properties
\begin{proposition}\label{rootspace-prop}(see \cite{Kn})

\vskip 0.2cm (i)
\;$\D\g=\g_{_0}\oplus\bigoplus_{\lambda\in\Sigma}\g_{_\lambda}$.

\vskip 0.2cm (ii)
\;$\D[\g_{_\lambda},\g_{_\gamma}]\subset\g_{_{\lambda+\gamma}}$
and if $\D\lambda\neq\gamma$, then
$\g_{_\lambda}\perp\g_{_\gamma}$ in the sense of $B$.

\vskip 0.2cm (iii) \;$\D\theta\g_{_\lambda}=\g_{_{-\lambda}}$, so
if $\D\lambda\in\Sigma$, then $-\lambda\in\Sigma$.

\vskip 0.2cm (iv) \;$\D\g_{_0}=\ba\oplus\am$ orthogonally.
\end{proposition}

\vskip 0.3cm Let $A=\exp\ba$. For any reductive Lie group $G$, it
has the following decomposition
\begin{theorem}\label{KAK-decomposition}(KAK decomposition \cite{Kn})
Every element in G has a decomposition as $k_{_1}ak_{_2}$ with
$k_{_1},k_{_2}\in K$ and $a\in A$. In this decomposition, $a$ is
uniquely determined up to conjugation by a member of $W$. If $a$
is fixed as $\exp H$ with $H\in \ba$ and $\lambda(H)\neq 0$ for
all $\lambda\in \Sigma$, then $k_{_1}$ is unique up to right
multiplication by a member of $Z_K(\ba)$.
\end{theorem}

\vskip 0.3cm
\subsection{Main result}

\quad\, Now let us present the main result of the paper. Let
$(G,K,\theta,B)$ be any reductive Lie group. Let $dg$ and $da$ be
the left-invariant measures on $G$ and $A$ respectively
corresponding to the Riemannian structure induced by $B$ and
$d\mu$ a $K\times K$-invariant measure on
 $(K\times K)/M$, where the set $M$ are
 defined by \eqref{M}.

\vskip 0.3cm
\begin{theorem}\label{weyl-int-mod} $\Big(K\times K,
\varphi, (G,dg), (A, da), \big((K\times K)/M, d\mu\big)\Big)$ is a
Weyl integration model with $d$-sheeted multiplicities, where
$d=|W|$ and the map $\varphi$ is defined by
 \eqref{reduced-map-varphi}.
\end{theorem}

\vskip 0.3cm
\begin{theorem}\label{main-result1} For any $f\in L^1(G,dg)$,
\begin{gather}\label{integral-formula1}
\int_Gf(g)\,dg=\int_A\Big(\int_{(K\times
K)/M}f(k_{_1}ak_{_2}^{-1})\,d\mu\Big) J(a)\,da,
\end{gather}
where $J(a)$ has the formula \eqref{noncompactJa}.
\end{theorem}

\vskip 0.5cm
\section{Proof of Main Results}

\vskip 0.3cm \quad\, In this section, we will prove the main
results by constructing a Weyl integration model and using the
theorem \ref{T:Weylbianhuan}.

\vskip 0.3cm Now Let $(G,K,\theta,B)$ be any reductive Lie group.
$\ba, A, \am, \Sigma, \Sigma^+$ and $W$ are defined in the above
section. We come to consider the G-action $\sigma$ defined by
\eqref{G-action-varphi} and its reduced map $\varphi$ defined by
\eqref{reduced-map-varphi}. Let $A'$ be the set of regular
elements in $A$ and $G'=KA'K$. It naturally has the following map
by restriction of $\varphi$, which is still denoted by $\varphi$.
\begin{gather}\label{reduced-map-varphi2}
\varphi:(K\times K)\big/M\times A'\rightarrow G'.
\end{gather}

\vskip 0.3cm
\begin{lemma}\label{M-ism-cenaK}
(i) $M$ is isomorphic to $Z_{K}(\ba)$.

(ii) $\varphi$ is a surjective $d$-sheeted map, where $d=|W|$.
\end{lemma}

\pf (i) Note that the unit element $e\in A$. If
$(k_{_1},k_{_2})\in M$, then $k_{_1}ek_{_2}^{-1}=e$, i.e.
$k_{_1}=k_{_2}$. Thus $M=\big\{(k,k)\in K\times K\,\big|\,ka=ak,
\forall\; a\in A\big\}$ is isomorphic to $Z_{K}(\ba)$.

(ii) By the conclusion $1$, it is the direct corollary of theorem
\ref{KAK-decomposition}.

\vskip 0.3cm Now for any $x\in G$, set
$M_x=\big\{(k_{_1},k_{_2})\in K\times
K\,\big|\,k_{_1}xk_{_2}^{-1}=x\big\}$ be the isotropic subgroup
associated with $x\in G$ and
$O_a=\big\{k_{_1}ak_{_2}^{-1}\,\big|\,(k_{_1},k_{_2})\in K\times
K\big\}$ the orbit of $a\in A'$. Then $M\subset M_a$, for all
$a\in A'$.

\vskip 0.3cm By the definition of reductive group, the invariant
bilinear $B$ determines an inner product on $\ap$, and $-B$
determines an inner product on $\ak$. We write $\ab=\ba^\perp$ in
$\ap$, and $\al=\am^\perp$ in $\ak$. By proposition
\ref{rootspace-prop}, it is obvious that
\begin{gather}\label{gablm}
\g=\am\oplus\al\oplus\ba\oplus\ab,\quad
\ab\oplus\al=\bigoplus_{\lambda\in\Sigma}\g_{_\lambda}, \quad
\al=\ak\bigcap\bigoplus_{\lambda\in\Sigma}\g_{_\lambda}, \quad
\ab=\ap\bigcap\bigoplus_{\lambda\in\Sigma}\g_{_\lambda}.
\end{gather}

For all $\lambda\in\Sigma^+$, we choose a normal orthogonal basis
$\{\xi_{_{\lambda,1}},\cdots,\xi_{_{\lambda,\beta_{_\lambda}}}\}$
in $\g_{_\lambda}$, where $\beta_{_\lambda}=\dim\g_{_\lambda}$
which may be larger than one. Then for all $\xi_{_{\lambda,j}}$,
we have
\begin{gather}
\theta(\xi_{_{\lambda,j}}+\theta \xi_{_{\lambda,j}})=
\xi_{_{\lambda,j}}+\theta \xi_{_{\lambda,j}}\in\al,\\
\theta(\xi_{_{\lambda,j}}-\theta
\xi_{_{\lambda,j}})=-(\xi_{_{\lambda,j}}-\theta
\xi_{_{\lambda,j}})\in\ab.
\end{gather}
Indeed, $\Bigl\{(\xi_{_{\lambda,j}}+\theta
\xi_{_{\lambda,j}})\,\Big|\,\lambda\in\Sigma^+,
 j=1,\cdots,\beta_{_\lambda}\Big\}$ composes of a basis of $\al$
and $\Bigl\{(\xi_{_{\lambda,j}}-\theta
\xi_{_{\lambda,j}})\,\Big|\,\lambda\in\Sigma^+,
j=1,\cdots,\beta_{_\lambda}\Bigr\}$ a basis of $\ab$. Then we get
\begin{gather}\label{dimalab}
\dim\al=\dim\ab=\sum_{\lambda\in\Sigma^+}\beta_{_\lambda}.
\end{gather}
Therefore,
\begin{align*}
\dim\big((K\times K)/M\times A\big)
&=\dim\big((\ak\oplus\ak)/\am\oplus\ba\big)
=\dim\ak+\dim\al+\dim\ba\\
&=\dim\ak+\dim\ab+\dim\ba=\dim\ak+\dim\ap=\dim G.
\end{align*}
So $\varphi$ is a map between the same dimension manifolds.

\vskip 0.3cm
\begin{lemma}\label{TaOa} (i) For $a\in A'$, denoted by $L_a$ the
left translation, then
\begin{gather}\label{TaOa2}
T_aO_a=\Big\{dL_a\big(Ad(a^{-1})\zeta_{_1}-\zeta_{_2}\big)\,\big|\,
(\zeta_{_1},\,\zeta_{_2})\in(\ak,\,\ak)\Big\}.
\end{gather}

(ii)  The following set is composed of a basis of $T_aO_a$,
\begin{gather}\label{basis-of-TaOa}
\begin{split}
F\equiv\Big\{dL_a(\eta_{_i}),\; dL_a(\xi_{_{\lambda,\,j}}^+),\;
dL_a(\xi_{_{\lambda,\,j}}^-)\;\big|\; i=1,\,2,\,\cdots,\,m,\\
\lambda\in\Sigma^+,\;j=1,\,2,\,\cdots,\,\beta_{_\lambda}\Big\}.
\end{split}
\end{gather}
\end{lemma}

\pf (i) Note that $T_aO_a$ is exactly the set composed of those
tangent vectors of the smooth curves $\exp(t\zeta_{_1})\cdot
a\cdot \exp(-t\zeta_{_2})$ at $t=0$ by the definition of orbit
$O_a$, where $\D (\zeta_{_1},\,\zeta_{_2})\in(\ak,\,\ak)$ and
$t<|\varepsilon|$. But
\begin{align}\label{tangent}
\begin{split}
\frac{d}{dt}\Big|_{t=0}\exp(t\zeta_{_1})\cdot a\cdot
\exp(-t\zeta_{_2})&=\frac{d}{dt}\Big|_{t=0}a\cdot
\exp\Big(t\big(Ad(a^{-1})\zeta_{_1}-\zeta_{_2}\big)+o(t^2)\Big)\\
&=dL_a\big(Ad(a^{-1})\zeta_{_1}-\zeta_{_2}\big).
\end{split}
\end{align}
Hence the conclusion is obtained.

 \vskip 0.3cm
(ii) Now suppose that
$\{\eta_{_1},\,\eta_{_2},\,\cdots,\,\eta_{_m}\}$ is a basis of
$\am$ and denote
$\xi_{_{\lambda,\,j}}^{\pm}\equiv\xi_{_{\lambda,j}}\pm\theta
\xi_{_{\lambda,j}}$. By $\ak=\am\oplus\al$, the following set is
composed of a basis of $(\ak,\,\ak)$,
\begin{gather}\label{basis-of-Tkkm}
\begin{split}
\Big\{(\eta_{_i},\,0),\;(\xi_{_{\lambda,\,j}}^+,\,0),\;
(0,\,\eta_{_i}),\;(0,\,\xi_{_{\lambda,\,j}}^+)\;\big|
\;i=1,\,2,\,\cdots,\,m,\;
\lambda\in\Sigma^+,\;j=1,\,2,\,\cdots,\,\beta_{_\lambda}\Big\}.
\end{split}
\end{gather}

Set $a=\exp H,\, H\in \ba$. Then by $[H,\,\eta_{_i}]=0$,
$[H,\,\xi_{_{\lambda,\,j}}^{\pm}]=\lambda(H)\xi_{_{\lambda,\,j}}^{\mp}$
and \eqref{tangent}, it can be calculated that
\begin{align}\label{Ada}
Ad(a^{-1})(\eta_{_i})-0&=\eta_{_i},\\
Ad(a^{-1})(\xi_{_{\lambda,\,j}}^+)-0
&=-\sinh\big(\lambda(H)\big)\xi_{_{\lambda,\,j}}^{-}
+\cosh\big(\lambda(H)\big)\xi_{_{\lambda,\,j}}^{+},\\
Ad(a^{-1})(0)-\eta_{_i}&=-\eta_{_i},\\
Ad(a^{-1})(0)-\xi_{_{\lambda,\,j}}^+&= -\xi_{_{\lambda,\,j}}^+.
\end{align}

Thus by the first conclusion (i), the set $F$ is exactly composed
of a basis of $T_aO_a$.

 \vskip 0.3cm
\begin{lemma}\label{dimM-Ma}
(i) $\dim M_a=\dim M$, for all $a\in A'$.

(ii) $T_aG=T_aA\oplus T_aO_a$, for all $a\in A'$.
\end{lemma}

\pf (i)  Since $M\cong Z_{_K}(\ba)$ in lemma \ref{M-ism-cenaK},
the Lie algebra of $M$ is $\am=Z_{\ak}(\ba)$. By lemma \ref{TaOa},
it is obvious that
\begin{gather}\label{TaOacong}
T_aO_a=dL_a(\ak\oplus\ab).
\end{gather}
Then by \eqref{dimalab},
\begin{align*}
\dim M_a&=\dim(K\times K)-\dim O_a\\
&=\dim\ak+\dim\ak-(\dim\ak+\dim\ab)\\
&=\dim\ak-\dim\al\\
&=\dim \am \\
&=\dim M.
\end{align*}

(ii) By \eqref{TaOacong}, it is directly obtained that
$$
T_aG=dL_a(\g)=dL_a(\ak\oplus\ap)
=dL_a\big(\ba\oplus(\ab\oplus\ak)\big)=T_aA\oplus
T_aO_a.
$$

\vskip 0.3cm \noindent{\bf Proof of theorem \ref{weyl-int-mod}.}
By Lemma \ref{M-ism-cenaK}-\ref{dimM-Ma}, proposition
\ref{P:localdiffeo} and the definition of Weyl integration model,
it is obvious.

\vskip 0.3cm \noindent{\bf Proof of theorem \ref{main-result1}.}
By theorem \ref{T:Weylbianhuan} and \ref{weyl-int-mod}, it only
needs to calculate $J(a)$ in \eqref{noncompact-Ja}. It is obvious
that $\D T_aA\perp T_aO_a$ for any $a\in A'$. Then by the second
conclusion in lemma \ref{dimM-Ma} and \eqref{P:det}, we need to
consider the map
\begin{gather}\label{Psia}
\Psi_a : T_{[(e,e)]}\big((K\times K)\big/M\big) \rightarrow
T_aO_a,\quad
(\zeta_{_1},\,\zeta_{_2})\mapsto\frac{d}{dt}\Big|_{t=0} \exp
(t\zeta_{_1})\cdot a\cdot\exp(-t\zeta_{_2}).
\end{gather}

Note that the tangent space $T_{[(e,e)]}\big((K\times
K)\big/M\big)$ of $(K\times K)\big/M$ at its unit element
$[(e,e)]$ is exactly isomorphic to
\begin{align}
(\ak,\al)&\equiv\Big\{(\zeta_{_1},\zeta_{_2})\,\big|\,
\zeta_{_1}\in\ak,\;\zeta_{_2}\in\al\Big\}.
\end{align}
Then the set
\begin{gather}\label{basis-of-Tkkm}
\begin{split}
\Big\{(\eta_{_i},\,0),\;(\xi_{_{\lambda,\,j}}^+,\,0),\;
(0,\,\xi_{_{\lambda,\,j}}^+)\;\big|\; i=1,\,2,\,\cdots,\,m,\;
\lambda\in\Sigma^+,\;j=1,\,2,\,\cdots,\,\beta_{_\lambda}\Big\}
\end{split}
\end{gather}
is composed of a basis of $T_{[(e,e)]}\big((K\times
K)\big/M\big)$.

\vskip 0.3cm It is completely analogous to the calculation in
\eqref{tangent} and \eqref{Ada}-(27), then we obtain
\begin{align}
\Psi_a\big((\eta_{_i},\,0)\big)
&=dL_{a}(\eta_{_i}),\\
\Psi_a\big((\xi_{_{\lambda,\,j}}^+,\,0)\big)
&=dL_{a}\big(-\sinh\big(\lambda(H)\big)\xi_{_{\lambda,\,j}}^{-}+
\cosh\big(\lambda(H)\big)\xi_{_{\lambda,\,j}}^{+}\big),\\
\Psi_a\big((0,\,\xi_{_{\lambda,\,j}}^+)\big)
&=dL_{a}(-\xi_{_{\lambda,\,j}}^+).
\end{align}
Thus by the second conclusion in lemma \ref{TaOa} and the
invariance of those chosen measures, there is a constant $C$ such
that
\begin{gather*}\label{noncompact-Ja}
J(a)=C\prod_{\lambda\in
\Sigma^+}\big|\sinh\big(\lambda(H)\big)\big|^{\beta_{_\lambda}}.
\end{gather*}
It completes the proof of theorem \ref{main-result1}.


\vskip 0.5cm

\begin{thebibliography}{9}


\bibitem{Kn} KNAPP, A. W., \textit{Lie Groups---Beyond
an Introduction}, 2nd edition,  Boston, Birkh\"auser,2002.

\bibitem{He2} HELGASON, S., \textit{Groups and Geometric
Analysis---Integral Geometry, Invariant Differetial Operators and
Spherical Functions}, Mathematical Surveys and Monographs, Vol.
83, New York, AMS, 2000.

\bibitem{AW} An, J. P. and Wang, Z. D., \textit{A Generalization
of Weyl Integration Formula}, Institute of Mathematics, Peking
University, Research Report No.18, 2004.

\bibitem{He} HELGASON, S., \textit{Differential Geometry,
Lie Groups, and Symmetric Spaces}, New tork, Acadymic Press, 1978.

\end{thebibliography}
\end{document}